\documentclass{amsart}
\usepackage{amsmath}
\usepackage{amsfonts}
\usepackage{amssymb}
\usepackage{amscd}

\DeclareMathOperator{\Z}{\mathbb Z}
\DeclareMathOperator{\B}{\mathcal B}
\DeclareMathOperator{\Sx}{\mathcal S}
\DeclareMathOperator{\Px}{\mathcal P}
\DeclareMathOperator{\Ob}{\mathcal Ob}
\DeclareMathOperator{\Q}{\mathbb Q}
\DeclareMathOperator{\C}{\mathbb C}
\DeclareMathOperator{\GL}{\bf GL}
\DeclareMathOperator{\Gm}{\bf G_m}
\DeclareMathOperator{\Ga}{\bf G_a}
\DeclareMathOperator{\Mn}{\bf M}
\DeclareMathOperator{\k0}{\bf k_0}
\DeclareMathOperator{\K}{\bf k}

\DeclareMathOperator{\diag}{diag}

\DeclareMathOperator{\Hom}{Hom}
\DeclareMathOperator{\Span}{span}
\DeclareMathOperator{\End}{\bf End}
\DeclareMathOperator{\id}{id}
\DeclareMathOperator{\ev}{ev}

\DeclareMathOperator{\CoDiffT}{\bf CoDiff}
\DeclareMathOperator{\Comod}{\bf Comod}

\DeclareMathOperator{\Aut}{\bf Aut}

\DeclareMathOperator{\IntEnd}{\underline{End}}

\DeclareMathOperator{\Rep}{\bf Rep}
\DeclareMathOperator{\Vect}{\bf Vect}
\DeclareMathOperator{\Cat}{\mathcal{C}}

\DeclareMathOperator{\Char}{char}

\DeclareMathOperator{\AlgT}{{\bf Alg}_{\k0}(\partial)}
\DeclareMathOperator{\Seq}{\mathcal{V}}
\DeclareMathOperator{\MSeq}{\mathcal{M}}

\theoremstyle{plain}
\newtheorem{theorem}{Theorem}
\newtheorem{lemma}{Lemma}
\newtheorem{proposition}{Proposition}
\newtheorem{corollary}{Corollary}

\theoremstyle{definition}
\newtheorem{definition}{Definition}
\newtheorem{example}{Example}

\theoremstyle{remark}
\newtheorem{remark}{Remark}

\newcommand{\Le}{\leqslant}
\newcommand{\Ge}{\geqslant}

\title[Tannakian categories and linear differential algebraic groups]{Tannakian categories, linear differential algebraic groups,
and parameterized linear differential equations}
\date\today
\author{Alexey Ovchinnikov}
\address{North Carolina State University, Department
of Mathematics, Box 8205, Raleigh, NC 27695-8205, USA}
\subjclass{12H05, 57T05, 18E99} \keywords{Tannakian categories,
linear differential equations,  Galois groups}
\urladdr{http://www.math.uic.edu/\textasciitilde aiovchin/}
\email{aiovchin@math.uic.edu}
\thanks{The work was partially supported by NSF Grant CCR-0096842
and by the Russian Foundation for Basic Research, project no.
05-01-00671.} \curraddr{University of Illinois at Chicago,
Department of Mathematics, Statistics, and Computer Science, 851 S.
Morgan Street, M/C 249, Chicago, IL 60607-7045, USA}
\begin{document}

\begin{abstract}
We provide conditions for a category with a fiber functor to be
equivalent to the category of representations of a linear differential
algebraic
group. This generalizes the notion of a neutral Tannakian category used
to
characterize the category of representations of a linear algebraic
group
\cite{Saavedra,Deligne}.
\end{abstract}

\maketitle

\section{Introduction}
Tannaka's Theorem (c.f., \cite{Springer}) states that a linear
algebraic
group is determined by its category of representations.  The problem of
recognizing when a category is the category of representations of a
linear
algebraic group (or more generally, an affine group scheme) is
attacked
via the theory of neutral Tannakian categories (see \cite{Saavedra,
Deligne}).  This theory allows one to detect the underlying presence of
a
linear algebraic group in various settings.  For example, the Galois
theory
of linear differential equations can be developed in this context (see
\cite{Deligne, DeligneFS, Michael}).

In   \cite{Cassidy}, Cassidy introduced the concept of a differential
algebraic group 
and in \cite{CassidyRep} studied the representation theory of linear
differential algebraic groups. Building on this work, we proved an
analogue
of Tannaka's Theorem for linear differential algebraic groups (see \cite{OvchRecoverGroup}). In the
present paper, we develop the notion of a neutral differential
Tannakian
category and show that this plays the same role for linear differential
algebraic groups that the theory of neutral Tannakian categories plays
for
linear algebraic groups.  As an application, we are able to give a
categorical development of the theory of parameterized linear
differential
equations that was introduced in \cite{PhyllisMichael}. 

Another approach to the Galois theory of systems
of linear differential equations with parameters is given in~\cite{TG}, where the authors
study Galois groups for generic values of the parameters. Also, it is shown in \cite{Hrushovsky} that over the field $\C(x)$ of rational functions over the complex numbers the differential Galois group of  a parameterized 
system of differential equations will be the same for all values of the parameter outside a countable 
union of proper algebraic sets.

The way we define differential Tannakian categories here relies on existence
of a fiber functor that ``commutes'' with the differential structure. It turns out
that there is a treatment of differential Tannakian categories not based on
the fiber functor \cite{DiffTann}, where there given conditions for a rigid abelian tensor category with an additional
differential structure so that the category has a fiber functor compatible with the differential structure. This is done based on the ideas of this paper, in the spirit of \cite{DeligneFS}, and extending \cite{Moshe}.

The paper is organized as follows. In Sections~\ref{Basics}
and~\ref{DefinitionProperties} we review the basic properties of linear
differential algebraic groups and define
and develop basic properties of a category which is a differential
analogue
of a neutral Tannakian category. In Section~\ref{MainSection} we prove
the main result of the paper: namely, there is a (pro-)linear differential
algebraic group such that a given neutral differential tannakian
category is the category of representations of the group. We give
an application of this result to the theory of parameterized linear
differential equations in Section~\ref{Applications}.
The techniques we use are built on
the
techniques introduced in \cite{Saavedra, Deligne} but the presence of
the
differential structure introduces new subtleties and several new
constructions.

\section{Basic Definitions}\label{Basics}
A $\Delta$-ring $R$,
where $\Delta = \{\partial_1,\ldots,\partial_m\}$, is a commutative associative ring with unit $1$ and commuting derivations $\partial_i : R\to R$ such that
$$
\partial_i(a+b) = \partial_i(a)+\partial_i(b),\quad \partial_i(ab) =
\partial_i(a)b + a\partial_i(b)
$$
for all $a, b \in R$. If $\k0$ is a field and
a $\Delta$-ring then $\k0$ is called a $\Delta$-field. We restrict ourselves to the case of
$$
\Char\k0 = 0.
$$ If $\Delta = \{\partial\}$ then
we call a $\Delta$-field as $\partial$-field.
For example, $\Q$ is a $\partial$-field with the unique
possible derivation (which is the zero one). The field
$\C(t)$ is also a $\partial$-field with $\partial(t) = f,$
and this $f$ can be any rational function in $\C(t).$ Let $C$ be
the field of constants of $\k0$, that is, $C = \ker\partial$. 

In Section~\ref{Applications} we require that every
consistent system of algebraic differential equations (that is,
it has a solution in an extension of $\k0$) with coefficients
in $\k0$ has a solution in $\k0$. Such a field is called differentially closed.
In characteristic zero for
a $\Delta$-field one can construct its differential algebraic
closure unique up to an isomorphism (see \cite[Definition 3.2]{PhyllisMichael} and the references given there). Also, in Section~\ref{Applications} we will deal with a field $\K$
equipped with two commuting differentiations, $\partial_x$
and $\partial_t$, and $\k0$ is the $\partial_t$-field of $\partial_x$-constants of $\K$ and
the differentiation $\partial_t$ that we use 
in Section~\ref{Applications} plays the role of
$\partial$ that we have here.

Let 
$$
\Theta = \left\{\partial^i\:|\: i\in \Z_{\Ge 0}\right\}.
$$
Since $\partial$ acts on a $\partial$-ring $R$,
there is a natural action of $\Theta$ on $R$.
A non-commutative ring $R[\partial]$ of linear differential operators is generated as a left $R$-module by the monoid $\Theta$. A typical element
of $R[\partial]$ is a polynomial 
$$
D = \sum_{i=1}^na_i\partial^{i},\ a_i \in R.
$$
The right $R$-module structure follows from the
formula 
$$
\partial\cdot a = a\cdot \partial + \partial(a)
$$
for all $a \in R$.
We denote the set of operators in $R[\partial]$ of order less than
or equal to $p$ by $R[\partial]_{\Le p}.$

Let $R$ be a $\partial$-ring. If $B$ is an $R$-algebra, then $B$ is a $\partial$-$R$-algebra
if the action of $\partial$ on $B$ extends the
action of $\partial$ on $R$. If $R_1$ and $R_2$
are $\partial$-rings then a ring homomorphism
$\varphi: R_1 \to R_2$ is called a $\partial$-homomorphism if it commutes with $\partial$, that is,
$$
\varphi\circ\partial = \partial\circ\varphi.
$$ We denote these homomorphisms simply by $\Hom(R_1,R_2)$.
If $A_1$ and $A_2$ are $\partial$-$\k0$-algebras
then a $\partial$-$\k0$-homomorhism simply means a 
$\k0[\partial]$-homomorphism. We denote the
category of $\partial$-$\k0$-algebras by $\AlgT$.
Let $Y = \{y_1,\ldots,y_n\}$ be a set of variables. We differentiate them:
$$
\Theta Y := \left\{\partial^iy_j
\:\big|\: i \in \mathbb{Z}_{\Ge 0},\ 1\Le j\Le n\right\}.
$$
The ring of differential polynomials $R\{Y\}$ in
differential indeterminates $Y$ 
over a $\partial$-ring $R$ is
the ring of commutative polynomials $R[\Theta Y]$
in infinitely many algebraically independent variables $\Theta Y$ with
the differentiation $\partial$, which naturally
extends $\partial$-action on $R$ as follows:
$$
\partial\left(\partial^i y_j\right) := \partial^{i+1}y_j
$$
for all $i \in \Z_{\Ge 0}$ and $1 \Le j \Le n$.
A $\partial$-$\k0$-algebra $A$ is called finitely
$\partial$-generated over $\k0$ if there exists
a finite subset $X = \{x_1,\ldots,x_n\} \subset A$
such that $A$ is a $\k0$-algebra generated by
$\Theta X$.  

An ideal $I$ in a $\partial$-ring $R$ is called differential if it is stable under the action of
$\partial$, that is,
$$
\partial(a) \in I
$$
for all $a \in I$. If $F \subset R$ then $[F]$ denotes the differential ideal generated by $F$.
If a differential ideal is radical, it is called
radical differential ideal. The radical differential
ideal generated by $F$ is denoted by $\{F\}$. If a
differential ideal is prime, it is called a prime
differential ideal.

\begin{definition}\label{SeqCategory}
The category $\Seq$ over a
$\partial$-field $\k0$ is the category of
finite dimensional vector spaces over $\k0$:
\begin{enumerate}
\item objects are finite dimensional $\k0$-vector spaces, 
\item morphisms are $\k0$-linear maps;
\end{enumerate}
with tensor product
$\otimes,$ direct sum $\oplus,$ dual $*,$ and additional operations:
$$
F^p : V \mapsto V^{(p)} := \k0[\partial]_{\Le p}\otimes V,
$$
which we call differentiation (or prolongation) functors. If $\varphi \in \Hom(V,W)$ then we define
$$
F^p(\varphi) : V^{(p)}\to W^{(p)},\ \varphi(\partial^q\otimes v) = \partial^q\otimes\varphi(v),\ 0\Le q\Le p.
$$
Here, $\k0[\partial]_{\Le p}$ is considered as the right $\k0$-module of differential operators up to order $p$ and $V$ is viewed as a left
$\k0$-module. We denote $F^1$ simply by $F$.
\end{definition}

For each $V \in \Ob(\Seq)$ there are: a natural inclusion
\begin{align}\label{eqincl}
i: V \to V^{(1)},\quad v \mapsto 1\otimes v,\:v \in V,
\end{align} a ``derivation''
\begin{align}\label{eqdiff}
\partial : V \to V^{(1)},\quad v \mapsto \partial\otimes v,
\end{align}
and a projection
\begin{align}\label{eqproj}\varphi : V^{(1)} \to V,\quad 1\otimes v \mapsto 0,\ \partial\otimes v \mapsto v.\end{align}

\begin{remark}
If $\{v_1,\ldots,v_n\}$ is a basis of $V$ then
$
\left\{v_1,\ldots,v_n,\ldots,\partial^p\otimes v_1,\ldots,\partial^p\otimes v_n\right\}
$
is a basis of $V^{(p)}.$
\end{remark}

\section{General definition of a differential Tannakian category}\label{DefinitionProperties}

\subsection{Definition}

\begin{definition}\label{NeutralParametricTannakian} A neutral differential Tannakian category $\Cat$ over a $\partial$-field  $\k0$ of characteristic zero is a
\begin{enumerate}
\item rigid
\item abelian
\item tensor
\end{enumerate}
category such that $\End(\underline{1})$
is the field $\k0$ supplied with an
\begin{enumerate}
\item exact
\item faithful
\item $\k0$-linear
\item tensor
\end{enumerate}
fiber functor $\omega: \Cat \to \Seq$ (see Definition~\ref{SeqCategory}); with a functor
$$
F : X \mapsto X^{(1)},
$$
intertwining\footnote{This means that there is a functorial in $X$ isomorphism $\omega\left(X^{(1)}\right) \cong \omega(X)^{(1)}$.} with $\omega$:
$$
F_{\Seq}\circ\omega = \omega\circ F_{\Cat},
$$
where the subscripts $\Cat$ and $\Seq$ indicate on which category
the functor $F$ acts\footnote{In Section~\ref{Applications} we will make use of $F^p$ on $\Seq$ for $p > 1$ for computational purposes.}.
 Also, for each $X \in \Ob(\Cat)$ there must be
an inclusion (a morphism with the trivial kernel)
$$
i: X \to X^{(1)}
$$
and surjection
$$
\phi: X^{(1)} \to X
$$
morphisms in $\Cat$,
so that the sequence
\begin{equation*}
\begin{CD}
0@>>> X @> i >> X^{(1)}@> \phi >> X @>>> 0
\end{CD}
\end{equation*}
is exact and
 is mapped by the fiber functor $\omega$ to the corresponding
exact sequence
in the category $\Seq$:
\begin{equation}\label{diffexact}
\begin{CD}
0 @>>> \omega(X) @>i>> \omega(X)^{(1)} @>\varphi>>\omega(X)@>>> 0
\end{CD},
\end{equation}
where the $\k0$-linear maps $i$ and $\varphi$ are defined in~\eqref{eqincl} and~\eqref{eqproj}, respectively.
Moreover, we require that the differential structure respects tensor
products (``product rule''), that is, for the natural $\k0$-linear map $$(\omega(X)\otimes\omega(Y))^{(1)} \to \omega(X)^{(1)}\otimes\omega(Y)^{(1)}$$
mapping
\begin{align*}
1\otimes (v\otimes u) &\mapsto (1\otimes v)\otimes(1\otimes u),\\
\partial\otimes (v\otimes u) &\mapsto (\partial\otimes v)\otimes (1\otimes u)+ (1\otimes v)\otimes(\partial\otimes u)
\end{align*}
there exists a corresponding morphism
in the category $\Cat$
\begin{equation}\label{difftensor}
\begin{CD}
(X\otimes Y)^{(1)}\\
@VV{}V\\
X^{(1)}\otimes Y^{(1)}
\end{CD}
\xrightarrow{\quad\omega\quad}
\begin{CD}
(\omega(X)\otimes\omega(Y))^{(1)}\\
 @VV{}V\\
 \omega(X)^{(1)}\otimes\omega(Y)^{(1)}
\end{CD}
\end{equation}
Also, we require that 
there exists
a morphism $g_X: (X^*)^{(1)}\to {X^{(1)}}^*$ satisfying
the following commutative diagram:
\begin{equation}\label{diffdual}
\begin{CD}
X^*@>i_{X^*}>>(X^*)^{(1)}@>\phi_{X^*}>>{X^*}  \\
@VV{\id}V@VV{g_X}V@VV{\id}V \\
X^*@>{\phi_X}^*>>\left(X^{(1)}\right)^*@>{i_X}^*>>{X^*}
\end{CD}
\end{equation}
such that the $\k0$-linear isomorphism $\omega(g_X) : \omega(X^*)^{(1)} \to \omega\left(X^{(1)}\right)^*$ is of the form:
\begin{equation}\label{diffdualomega}
\begin{CD}
\omega(X)^*@>i_{\omega(X)^*}>>(\omega(X)^*)^{(1)}\cong\omega(X)^*\oplus\omega(X)^*@>\varphi_{\omega(X)^*}>>{\omega(X)^*}  \\
@VV{\id}V@VV{\omega(g_X)}V@VV{\id}V \\
\omega(X)^*@>{\varphi_{\omega(X)}^*}>>\left(\omega(X)^{(1)}\right)^*\cong(\omega(X)\oplus\omega(X))^*@>{i_{\omega(X)}^*}>>{\omega(X)^*}
\end{CD}
\end{equation}
with the map $\omega(g_X)$ respecting the 
given splittings into direct sums. 
\end{definition}

In Example~\ref{RepExample}, Section~\ref{secdiffcomod}, we will show that the category of
differential representations of a linear differential algebraic
group together with the forgetful functor form a neutral differential Tannakian category.

\begin{remark}
The functor $F$ and these additional morphisms are used directly to recover
the differential structure on the (pro-)linear differential algebraic group
whose category of representations is $\Cat.$ This explains why
they are given here. In particular, conditions~\eqref{diffdual} and \eqref{diffdualomega} are given to ensure that the coinverse on the Hopf algebra that we recover is a {\it differential} homomorphism. This will further be discussed in Lemma~\ref{ComputationalDiffDual} and Proposition~\ref{FinitelyGenerated}.
\end{remark}
Our main goal is to prove the following theorem.

\begin{theorem}\label{FirstFormulation} For a neutral differential Tannakian category $\Cat$:
\begin{itemize}
\item the functor $\Aut^{\otimes,\partial}(\omega)$ from $\AlgT$ to $\{\mathrm{Groups}\}$ is representable by a differential Hopf algebra $A,$
\item the functor $\omega$ defines an equivalence of tensor categories
$\Cat \to \Rep_G,$ with $G$ being the affine differential group
scheme represented by $A.$
\end{itemize}
\end{theorem}

We start with developing a technique for this. We partially follow \cite[pages 130--137]{Deligne} modifying their definitions and proofs to give the correct result in the differential case. 
In Example~\ref{ExampleFromRepToAut}, Section~\ref{secdiffstructure}, we will show the ``easy part'' of Theorem~\ref{FirstFormulation} when $\Cat = \Rep_G$.

\subsection{Tensor product $\Cat\otimes \Seq$ and its properties}
In this section we will
list several lemmas given in \cite[pages 131--132]{Deligne} that we need to use.
Let $(\Cat, \omega)$ be a neutral differential Tannakian category. Consider objects $V \in \Ob(\Seq)$ and $X \in \Ob(\Cat).$

\begin{definition}
We define
$$
V \otimes X = \{[(X^n)_\alpha, \phi_{\beta,\alpha}]\:|\:\alpha : (\k0^n) \cong V\},
$$
where $(X^n)_\alpha := X\oplus\ldots\oplus X$ and $\phi_{\beta,\alpha} : (X^n)_\alpha \to (X^n)_\beta$
is defined by $\beta^{-1}\circ\alpha$ acting
as an element of $\GL_n(\k0).$
\end{definition}

\begin{lemma} $\phi_{\gamma,\beta}\circ\phi_{\beta,\alpha} = \phi_{\gamma,\alpha}.$
\end{lemma}

\begin{definition} For an object $T \in \Ob(\Cat)$ we say that
$\Phi \in \Hom(V\otimes X, T)$ if $\Phi = \{\varphi_\alpha : (X^n)_\alpha \to T\}$ making possible
diagrams with $\phi_{\beta,\gamma}$ commutative. Also, $\Phi \in \Hom(T,V\otimes X)$
if $\Phi = \{(\varphi_1,\ldots,\varphi_n)_\alpha\},$ $\varphi_i : T \to X$ making possible
diagrams with $\phi_{\beta,\gamma}$ commutative.
\end{definition}

Recall that $\k0 \subset \End(X).$

\begin{lemma}\label{FromVtoHom} For $V \in \Ob(\Seq)$ there is a canonical $\k0$-linear map
\begin{align*}
V \to \Hom(X,&V\otimes X), \quad v \mapsto ((\varphi^1,\ldots,\varphi^n)_\alpha : X \to X_\alpha),\\
&(\varphi^i)_\alpha\ \text{is the multiplication by}\ w^i,
\end{align*}
where $\alpha^{-1}(v) = w,\ w = (w^1,\ldots,w^n)\in \k0^n.$
\end{lemma}

\begin{lemma}\label{HomIsom} We have a functorial isomorphism
$$
\Hom(V\otimes X,T) \cong \Hom\left(V, \Hom(X,T)\right),
$$
where $T$ is an object in $\Cat$.
\end{lemma}

\begin{lemma}\label{ApplyFunctor} Let $F : \Cat \to \Cat'$ be a functor. Then $F(V\otimes X) = V\otimes F(X)$.
\end{lemma}

\begin{lemma}\label{UsualOnSeq} Let $\Cat = \Seq$ then $V\otimes X \cong V\otimes_{\Seq}X$ (the
tensor product in $\Seq$).
\end{lemma}

\subsection{Differential comodules}\label{secdiffcomod}
 Let $A$ be a $\partial$-$\k0$-algebra.
Assume that $A$ is supplied with the following operations:
\begin{itemize}
\item differential algebra homomorphism $m : A\otimes A \to A$ is the
multiplication map on $A,$
\item differential algebra homomorphism $\Delta : A \to A\otimes A$
which is a comultiplication,
\item differential algebra homomorphism $\varepsilon : A  \to \k0$
which is a counit,
\item differential algebra homomorphism $S : A \to A$
which is a coinverse.
\end{itemize}
We also assume that these maps satisfy commutative diagrams (see \cite[page 225]{CassidyRep}):
\begin{equation}\label{HopfAxioms}
\begin{CD}
A @>{\Delta}>> A\otimes A\\
@VV{\Delta}V @VV{\id_A\otimes\Delta}V \\
A\otimes A @>{\Delta\otimes \id_A}>> A\otimes A\otimes A
\end{CD}\qquad
\begin{CD}
A @>{\Delta}>> A\otimes A\\
@VV{\id_B}V @VV{\id_A\otimes\varepsilon}V \\
A @>{\sim}>> A\otimes \k0
\end{CD}\quad\quad
\begin{CD}
A @>{\Delta}>> A\otimes A\\
@VV{\varepsilon}V @VV{m\circ(S\otimes \id_A)}V \\
\k0 @>{\hookrightarrow}>> A
\end{CD}
\end{equation}

\begin{definition}\label{DiffHopfAlgebra} Such a commutative associative $\partial$-$\k0$-algebra $A$
with unit $1$ and operations $m,$ $\Delta,$ $S,$
and $\varepsilon$ satisfying axioms~\eqref{HopfAxioms} is called a {\it differential Hopf algebra} (or $\partial$-$\k0$-Hopf algebra).
\end{definition}

\begin{definition}\label{DiffCoalgebra} Assume that $A$ is just a vector space over $\k0$ with a derivation $\partial$ which extends the
derivation on $\k0$ equipped with $\Delta$ 
and $\varepsilon$
commuting with $\partial.$ In this case $A$ is called a {\it differential coalgebra}. When, in addition, $A$ has $S,$
it is called a {\it differential bialgebra}.
\end{definition}

\begin{definition}\label{DiffComodules} A  finite dimensional vector space $V$ over $\k0$ is called a {\it $A$-comodule}
if there is a given $\k0$-linear morphism
$$
\rho : V \to V\otimes A,
$$
satisfying the axioms:
\begin{equation*}
\begin{CD}
V@>{\rho}>> V\otimes A\\
@VV{\rho}V @VV{\id_V\otimes\Delta}V \\
V\otimes A@>{\rho\otimes\id_A}>> V\otimes A\otimes A
\end{CD}\qquad\qquad
\begin{CD}
V @>{\rho}>> V\otimes A\\
@VV{\id_V}V @VV{\id_V\otimes\varepsilon}V \\
V @>{\sim}>> V\otimes \k0
\end{CD}
\end{equation*}
If $A$ is a differential coalgebra with $\Delta$ and $\varepsilon$
then for a comodule $V$ over $A$ the $\k0$ space $\k0[\partial]_{\Le i}\otimes V = V^{(i)}$ has a natural $A$-comodule structure
\begin{align}\label{DiffFla}
\rho(f\otimes v) := f\otimes\rho(v)
\end{align}
for  $f \in \k0[\partial]_{\Le i}$ and $v \in V.$
\end{definition}

We denote the category of  comodules over a differential coalgebra
$A$ by $\CoDiffT_A$ with the induced differentiation~\eqref{DiffFla}.
For a $\partial$-$\k0$-Hopf algebra $A$ the functor
$$
G : \AlgT \to \{\mathrm{Groups}\}, \quad R \mapsto \Hom(A,R)
$$
is called an {\it affine differential algebraic group scheme} generated by $A$ (see \cite[Section 3.3]{OvchRecoverGroup}). In this case $V \in \Ob(\CoDiffT_A)$ is called a
{\it differential representation} of $G$ (\cite[Definition 7, Theorem 1]{OvchRecoverGroup}).
The category $\CoDiffT_A$ is also denoted by $\Rep_G$.

The differential $\GL_n$ by definition is the functor 
represented by the $\partial$-$\k0$-Hopf algebra
$$
\k0\{X_{11},\ldots,X_{nn},1/\det(X)\},
$$
where $X_{ij}$ are differential indeterminates.
The comultiplication $\Delta$ and coinverse  $S$ are defined on $X_{ij}$ in the usual way. Their
prolongation on the derivatives of $X_{ij}$ can be
obtained by differentiation.
\begin{example} Denote the differential $\GL_1$ by 
$\Gm$. The $\partial$-$\k0$-Hopf algebra for $\Gm$
is then $\k0\{y,1/y\}$. We have:
\begin{align*}
\Delta(y) &= y\otimes y,\\
S(y) &= 1/y.
\end{align*}
These maps are $\partial$-homomorphisms. Therefore,
\begin{align*}
\Delta(\partial y) = \partial(\Delta(y))=\partial y\otimes y + y\otimes\partial y,\\
S(\partial y) = \partial(S(y)) = \partial(1/y) = -\partial y/y^2
\end{align*}
and so on.
\end{example}
Recall that $C$ is the field of $\partial$-constants of $\k0$.
\begin{example} The ``constant'' multiplicative group $\Gm(C)$ is the subgroup of $\Gm$ given
by the equation $\partial y=0$. So, the $\partial$-$\k0$-Hopf algebra of $\Gm(C)$ is
$$
\k0\{y,1/y\}/[\partial y] \cong \k0[y,1/y],
$$
which is the usual algebraic multiplicative group.
\end{example}

\begin{example} The differential $\Gm$ has another
non-trivial differential algebraic subgroup given by the following differential equation
$$
\partial\left(\frac{\partial y}{y}\right) = 0
$$
(see \cite[page 126]{PhyllisMichael} for more information).
\end{example}

\begin{example} The differential $\Ga$ is represented by the $\partial$-$\k0$-Hopf algebra $\k0\{y\}$ 
with 
\begin{align*}
\Delta(\partial^py) &= \partial^py\otimes 1 + 1\otimes\partial^py,\\
S(\partial^py) &= (-1)^{p+1}\partial^py
\end{align*}
for all $p \in\Z_{\Ge 0}$.
\end{example}

\begin{definition} An affine differential algebraic group scheme $G$ is called a {\it linear
differential algebraic group} if there exists
an imbedding:
$$
G \to \GL_n
$$
for some $n \in \Z_{\Ge 1}$.
\end{definition}

\begin{example}\label{RepExample} We show that for a linear differential algebraic group $G$ the category $\Rep_G$ together
with the forgetful functor form a neutral differential Tannakian
category. Indeed, in the above it is shown how to define the
functor $F$.
Moreover, for any $g \in G$, $V$ and $W \in \Ob(\Rep_G)$, $v\in V$,
$w\in W$ we have:
$$
V \ni g\cdot v \mapsto 1\otimes (g\cdot v) = g\cdot (1\otimes v) \in V^{(1)};
$$ in addition, applying $\partial : V \to
V^{(1)}$ we obtain
$$
g\cdot v \mapsto \partial\otimes (g\cdot v) = g\cdot(\partial\otimes v);
$$ 
also, for the product rule we have:
\begin{align*}
(V\otimes W)^{(1)} \ni\: & g\cdot(\partial\otimes (v\otimes w))
=\partial\otimes(g\cdot v\otimes g\cdot w) \mapsto \\
&\mapsto (\partial\otimes g\cdot v)\otimes (1\otimes g\cdot  w)+
(1\otimes g\cdot v)\otimes(\partial\otimes g\cdot w)=\\
&= g\cdot((\partial\otimes v)\otimes (1\otimes w)+(1\otimes
v)\otimes(\partial\otimes w))\in V^{(1)}\otimes W^{(1)}.
\end{align*}
Moreover,
$$
V^{(1)}\ni g\cdot \partial\otimes v = \partial\otimes (g\cdot
v)\mapsto g\cdot v \in V.
$$
Finally, let $V=\Span_{\k0}\{v_1,\ldots,v_n\}$, $V^* =
\Span_{\k0}\{v_1^*,\ldots,v_n^*\}$, $v = \sum a_i v_i$, and  $g\cdot
v_j^* = \sum g_{ji}v_i^*$. We then have
$$v_j^*\left(g^{-1}v\right) =: v_j^*\left(\sum c_k v_k\right) = c_j.$$
Hence,
$$c_j = v_j^*\left(g^{-1}v\right) = (g\cdot v_j^*)(v) = \sum g_{ji}v_i^*(v) = \sum g_{ji}a_i.$$
\begin{align*}
(V^*)^{(1)} \ni g\cdot (\partial\otimes v_j^*) &= \partial\otimes (g\cdot v_j^*) = \partial\otimes \sum g_{ji}v_i^* = \\
&=\sum\left(g_{ji}\partial\otimes v_i^*+ \partial(g_{ji})\otimes v_i^*\right)\mapsto \\
&\mapsto\sum\left(g_{ji}F(v_i^*)+\partial(g_{ji})\otimes (\partial\otimes v_i)^*\right) \in \left(V^{(1)}\right)^*.
\end{align*}
Now,
\begin{align*}
\sum\left(g_{ji}F(v_i^*)+\partial(g_{ji})\otimes (\partial\otimes v_i)^*\right)(1\otimes v) &= \sum g_{ji}a_i = c_j = \\
&=F(v_j^*)(g^{-1}\cdot (1\otimes v)) =\\
&=
\left(g\cdot F(v_j^*)\right)(1\otimes v)
\end{align*}
and
\begin{align*}
\sum\left(g_{ji}F(v_i^*)+\partial(g_{ji})\otimes (\partial\otimes v_i)^*\right)(\partial\otimes v) &= \sum g_{ji}\partial(a_i)+\partial(g_{ji})a_i =\\
&= \partial\left(\sum g_{ji}a_i\right) = \\
&=\partial c_j = \partial\left(v_j^*\left(g^{-1}\cdot v\right)\right) = \\
&=
F(v_j^*)\left(\partial\otimes\left(g^{-1}\cdot v\right)\right) = \\
&=(g\cdot F(v_j^*))(\partial\otimes v).
\end{align*}
In addition,
\begin{align*}
(V^*)^{(1)}\ni g\cdot (1\otimes v_j^*) = \sum g_{ji}\otimes v_i^*\mapsto \sum g_{ji} (\partial\otimes v_i)^*
\end{align*}
and
\begin{align*}
&&\sum g_{ji}(\partial\otimes v_i)^*(1\otimes v) &= 0 =
\sum(\partial\otimes v_j)^*(c_i\otimes v_i)=g\cdot(\partial\otimes v_j)^*(1\otimes v),\\
&&\sum g_{ji} (\partial\otimes v_i)^*(\partial\otimes v) &= \sum g_{ji} a_i = c_j = (\partial\otimes v_j)^*\left(\partial\otimes \left(g^{-1}\cdot v\right)\right)=\\
&&&=g\cdot(\partial\otimes v_j)^*(\partial\otimes v).
\end{align*}
This shows that all these maps are $G$-morphisms, which is all
that we needed.
\end{example}

\subsection{Recovering a bialgebra over $\k0$}

\begin{remark} Let $\Cat = \Rep_G$ for a linear differential
algebraic group $G$ and $\omega$ is the forgetful functor.
This is one of main applications of the theory we are developing.
It turns out that one can explicitly recover the group $G$ knowing $\Cat$.
The paper \cite{OvchRecoverGroup} does this.
\end{remark}

Let $X$ be an object of $\Cat$ and $\{\{X\}\}$ (respectively, $\Cat_X$) be the full abelian (respectively, full abelian tensor) subcategory of $\Cat$ generated by $X$
(respectively, containing $X$ and closed under $F$).

Denote by $B_X$ the object $P$  (or $P'$) contained in $X^*\otimes\omega(X)$ and constructed as
in \cite[Lemma 2.12]{Deligne} for the category $\{\{X\}\}$. Its image $\omega(B_X)$ stabilizes all $\omega(Y),$
where $Y$ is an object of $\{\{X\}\}.$  Denote $P_X = (B_X)^*.$ We note that $A_X := \omega(P_X)$ is a finite
dimensional vector space over $\k0.$
One can construct $P_X$ directly as follows.
Consider
$$
F_X = \bigoplus_{V \in \Ob(\{\{X\}\})}  V\otimes\omega(V)^*.
$$
For an object $V$ of $\{\{X\}\}$ we have the canonical injections:
$$
i_V : V\otimes\omega(V)^* \to F_X.
$$
Consider the minimal subobject $R_X$ of $F_X$ with subobjects
$$ \left\{\left(i_V(\id\otimes \phi^*) -
i_W(\phi\otimes \id)\right)(V\otimes
\omega(W)^*)\big|V,W \in\Ob(\{\{X\}\}),\:\phi\in
\Hom\left(V,W\right)\right\}.
$$
We let
$$
P_X = F_X\big/R_X.
$$
We now put
$$
A_X = \omega(P_X),
$$
which is a $\k0$-vector space (see Lemma~\ref{ApplyFunctor}).
Let $V$ be an object of $\{\{X\}\}.$
For
$$
v  \in \omega(V),\
u \in \omega(V)^*
$$
we denote by
$$
a_V\left(v\otimes u\right)
$$
the image in $A_X$ of
$$
\omega(i_V)\left(v\otimes u\right).
$$
So, for any $\phi \in \Hom(V,W)$ we have
\begin{align}\label{morphismsanda}
a_V\left(v\otimes \omega(\phi)^*(u)\right) =
a_W(\omega(\phi)(v)\otimes u).
\end{align}
for all
$$
v \in \omega(V),\ u \in \omega(W)^*.
$$

Let us define a comultiplication on $A_X.$
Let $\{v_i\}$ be a basis of $\omega(V)$ with the dual basis
$\{u_j\}$ of $\omega(V)^*.$
We let
\begin{equation}\label{Comult}
\Delta : a_V(v\otimes u)
\mapsto\sum a_V(v_i\otimes u)\otimes a_V(v\otimes u_i).
\end{equation}
The counit is defined in the following way:
\begin{equation}\label{Counit}
\varepsilon : a_V(v\otimes u) \mapsto u(v).
\end{equation}
The coinverse is defined as follows:
\begin{equation}\label{Coinverse}
S : a_V(v\otimes u) \mapsto a_{V^*}(u\otimes v).
\end{equation}

\begin{proposition}\label{bialgebra} With the operations~\eqref{Comult},  \eqref{Counit}, and \eqref{Coinverse} the finite dimensional $\k0$-vector space $A_X$ is a bialgebra.
\end{proposition}
\begin{proof}
This follows directly from \cite[Lemmas 9, 10, and
11]{OvchRecoverGroup}. In these lemmas it is shown that the
operations given by formulas~\eqref{Comult}, \eqref{Counit}, and
\eqref{Coinverse} are well defined and satisfy the usual axioms
of a bialgebra. Moreover, in the proofs only the fact that the
category is rigid and abelian was used.
\end{proof}

\subsection{Differential structure}\label{secdiffstructure}

\begin{definition}\label{TensorAutomorphism} For a $\partial$-$\k0$-algebra $R$
we define a group $\Aut^{\otimes,\partial}(\omega)(R)$ to be the set of sequences
$$\lambda(R) = (\lambda_X\:|\: X\in \Ob(\Cat)) \in
\Aut^{\otimes,\partial}(\omega)(R)$$ such that $\lambda_X$ is an
$R$-linear automorphism of $\omega(X)\otimes R$ for each object $X$,
$\omega(X) \in \Ob(\Seq)$, that is, $\lambda_X \in
\Aut_R(\omega(X)\otimes R)$, such that
\begin{itemize}
\item for all $X_1,$ $X_2$ we have
\begin{align}\label{TensorSpreading}
\lambda_{X_1\otimes X_2} = \lambda_{X_1}\otimes\lambda_{X_2},
\end{align}
\item $\lambda_{\underline{1}}$ is the identity map on
$\underline{1}\otimes R = R$,
\item for every $\alpha \in \Hom(X,Y)$ we have
\begin{align}\label{Equivariance}
\lambda_Y\circ(\alpha\otimes \id_R) = (\alpha\otimes
\id_R)\circ\lambda_X : \omega(X)\otimes R \to \omega(Y)\otimes R,
\end{align}
\item for every $X$ we have
\begin{align}\label{CommuteWithD}
\partial\circ\lambda_X = \lambda_{X^{(1)}}\circ \partial,
\end{align}
\item the group operation $\lambda_1(R)\cdot\lambda_2(R)$ is defined
by composition in each set $\Aut_R(\omega(X)\otimes R).$
\end{itemize}
\end{definition}

\begin{example}\label{ExampleFromRepToAut} For $\Cat = \Rep_G$,
where $G$ is a linear differential algebraic group, we will construct an injection $G \to \Aut^{\otimes,\partial}(\omega).$
For this take any differential $\k0$-algebra $R$ and $g \in G(R).$
For any $V \in \Rep_G$ the element $g$ defines an $R$-linear automorphism
$$\lambda_V : \omega(V)\otimes R \to \omega(V)\otimes R.$$
We show that $\{\lambda_V\}$ is an element of $\Aut^{\otimes,\partial}(\omega)(R).$ Indeed, formulas \eqref{TensorSpreading}, \eqref{Equivariance}, and injectivity are the usual properties
of representation. We show that formula \eqref{CommuteWithD} holds true as well:
\begin{align*}
\lambda_{X^{(1)}}\circ\partial (v\otimes r) =& \lambda_{X^{(1)}}((\partial\otimes v)\otimes r + (1\otimes v)\otimes\partial r) =
(\partial\otimes (g\cdot v)) \otimes r +\\
&+  (1\otimes (g\cdot v))\otimes\partial r) = \partial ((g\cdot v)\otimes r) = \partial\circ\lambda_X(v\otimes r)
\end{align*}
for any $r \in R$.
The less obvious part is surjectivity and instead of showing
it directly we prove the main result using the dual language of
the corresponding differential Hopf algebras.
\end{example}

Speaking informally, in the previous sections things looked quite
algebraic without differential algebra. The differential part
appears here: we have an additional condition~\eqref{CommuteWithD},
which allows one to prove the differential analogue of Tannaka's
theorem.
We are now going to look at the ``differential'' category $\Cat_X$
generated by $X$. For $X \in \Ob(\Cat)$ we let
$$\Px_X = \varinjlim_{Y\in\Ob(\Cat_X)}P_Y.$$
Let
$$
\B_X = \omega(\Px_X).
$$
We {\bf finally} put
$$
A = \varinjlim_{X \in \Ob(\Cat)}\B_X.
$$
This is a coalgebra due to the
previous considerations. We are going to make it a {\it differential} Hopf algebra using 
the prolongation functor $F$
and $\otimes$.

For this take
$v \in \omega(V),$ $w \in \omega(W),$ $u \in \omega(V)^*,$
$t \in \omega(W)^*$
and let
\begin{equation}\label{Mult}
a_V(v\otimes u)\cdot a_W( w\otimes t) = a_{V\otimes W}((v\otimes w))\otimes(u\otimes t)).
\end{equation}
Define
$$F(u)(1\otimes v) := u(v),\ F(u)(\partial\otimes v) :=
\partial(u(v)),\ u \in V^*,\: v\in V$$
for a vector space $V$. 
Introduce a
{\it differential} structure on $A:$
\begin{equation}\label{Diff}
\partial\left(a_V(v\otimes u)\right) = a_{V^{(1)}}((\partial v)\otimes F(u)),
\end{equation}
which can be extended to the whole $A$ as a derivation. Recall that
we extend here $u \in V^*$ to an element of $\left(V^{(1)}\right)^*$
by
$$F(u)(\partial\otimes w) := \partial(u(w)),\quad F(u)(1\otimes w) := u(w).$$
Note also that
$$ a_V(v\otimes u) = a_{V^{(1)}}\left(\partial
v\otimes \varphi^*(u)\right).
$$
We recall here that $\varphi^*$ is the dual map to the morphism
given in~\eqref{eqproj}:
$$\phi: V^{(1)} \to V,\quad \varphi := \omega(\phi) : 1\otimes v \mapsto 0,\ \partial\otimes v\mapsto v.$$

The following ``computational'' version of~\eqref{diffdual} and \eqref{diffdualomega}
will further be used in Proposition~\ref{FinitelyGenerated}. 

\begin{lemma}\label{ComputationalDiffDual} 
Let $V\in \Ob(\Cat_X)$ and $\{v_1,\ldots,v_n\}$
be an ordered basis of $\omega(V)$ corresponding to the splitting in~\eqref{diffdualomega} with the dual basis $\{v_1^*,\ldots,v_n^*\}$ of the dual vector space  $\omega(V)^*$. Then,
$$
a_{{(V^*)}^{(1)}}\left(\partial (v_j^*)\otimes F(v_i)\right) = a_{{\left(V^{(1)}\right)}^*}(F(v_j^*)\otimes\partial
v_i).
$$
\end{lemma}
\begin{proof}Let $g_X$ be the morphism $(V^*)^{(1)}\to {V^{(1)}}^*$ given by~\eqref{diffdual}.
We then have $$
i_{\omega(V)^*}(v_i^*)= 1\otimes v_i^*
$$ 
and
$${\varphi_{\omega(V)}}^*(v_i^*)(1\otimes v_j)= 0,\quad {\varphi_{\omega(V)}}^*(v_i^*)(\partial\otimes v_j)=\delta_{ij}.$$ Therefore,
\begin{align}\label{phidual}
{\varphi_{\omega(V)}}^*(v_i^*) = (\partial\otimes v_i)^*
\end{align}
and
$$
\omega(g_V)(1\otimes v_i^*) = {\varphi_{\omega(V)}}^*(v_i^*)= (\partial\otimes v_i)^*.
$$
Now,
$$
\varphi_{\omega(V)^*}(\partial\otimes v_i^*) = v_i^*
$$
and for any $a_1,\ldots,a_n \in \k0$ we have
$$
{i_{\omega(V)}}^* : (1\otimes v_i)^* + \sum_{j=1}^na_j\cdot(\partial\otimes v_j)^* \mapsto v_i^*.
$$
Hence, there exist $b_{i1},\ldots,b_{in} \in \k0$ such that
$$
\omega(g_V)(\partial\otimes v_i^*) = (1\otimes v_i)^* + \sum_{j=1}^nb_{ij}\cdot(\partial\otimes v_j)^* .
$$
Since $\omega(g_V)$ respects the splittings into direct sums given in~\eqref{diffdualomega}, $b_{ij} = 0$ for all $1\Le i,j\Le n$.
Note that
$$F(v_i^*)(1\otimes v_j) = \delta_{ij},\quad F(v_i^*)(\partial\otimes v_j) = \partial(v_i^*(v_j)) = 0.$$ Therefore,
$$
F(v_i^*) = (1\otimes v_i)^* = \omega(g_V)(\partial\otimes v_i^*).
$$
Moreover,
$$
\omega(g_V)^*(\partial v_i)(1\otimes v_j^*) = (\partial v_i)(\omega(g_V)(1\otimes v_j^*))=
\partial v_i((\partial v_j)^*) = \delta_{ij}
$$
and
$$
\omega(g_V)^*(\partial v_i)(\partial\otimes v_j^*) = (\partial v_i)(\omega(g_V)(\partial\otimes v_j^*))= 
(1\otimes v_j)^*(\partial v_i) = 0.
$$
Hence,
$$
\omega(g_V)^*(\partial v_i) = (1\otimes v_i^*)^* $$
Now,
$$
F(v_i)(1\otimes v_j^*)= \delta_{ij},\quad F(v_i)(\partial\otimes v_j^*) = 0.
$$
Therefore,
$$
F(v_i) = (1\otimes v_i^*)^* = \omega(g_V)^*(\partial v_i),
$$
which is what we finally needed to apply~\eqref{morphismsanda}.
\end{proof}

\begin{proposition}\label{FinitelyGenerated} With the operations~\eqref{Mult}, \eqref{Diff}, \eqref{Comult}, \eqref{Counit},
and \eqref{Coinverse} the $\k0$-vector space $A$ is a direct limit of finitely generated commutative associative differential Hopf algebras with the unity.
\end{proposition}
\begin{proof}
All the statements follow from Proposition~\ref{bialgebra} and
 \cite[Lemmas 5, 6, 7, 9, 10, and 11]{OvchRecoverGroup} except for finitely generated, because our construction of $\B_X$ and operations
on it are the same as in \cite{OvchRecoverGroup}.
The $\partial$-$\k0$-algebra $\B_X$ is generated by the
elements $a_X(v\otimes u).$  The statement now follows as
$\omega(X)\otimes\omega(X)^*$ is a finite dimensional $\k0$-vector space.

In order to show the other properties like:
\begin{itemize}
\item $\partial,\:m,\:\Delta,\:\varepsilon,\:S$ are well defined; 
\item commutation with $\partial$ of $\Delta,\:S,\:\varepsilon$;
\item respecting multiplication ($\Delta,\:S,\:\varepsilon$ are algebra homomorphisms);
\item the product rule for the multiplication,
\end{itemize}
in these lemmas only the properties of a neutral differential
Tannakian category were used. In particular, properties
\eqref{diffexact}, \eqref{difftensor},  \eqref{diffdual}, and \eqref{diffdualomega} were needed for the proofs. That is why one finds
them
among the axioms of a neutral differential Tannakian category. Moreover, we have Lemma~\ref{ComputationalDiffDual}, which is used in \cite[Lemma 11]{OvchRecoverGroup} to show that the coinverse we have defined commutes with differentiation $\partial$.
\end{proof}

\subsection{Equivalence of categories}

\begin{lemma}\label{FactorsThrough} The restriction $\omega|_{\Cat_X} : \Cat_X \to \Seq$
factors through $\CoDiffT_{\B_X}.$
\end{lemma}
\begin{proof}
Let $Y$ be an object of $\Cat_X.$ We introduce a $\B_X$-comodule
structure on $\omega(Y).$
Let $R$ be a $\partial$-$\k0$-algebra.
Consider $\xi \in \Hom(\B_X,R).$ Let $v \in \omega(Y)$ and $u \in\omega(Y)^*.$ There is an endomorphism $\lambda_Y$ of $\omega(Y)\otimes R$
such that
\begin{align}\label{FactorFormula}
\langle\lambda_Y(v),u\rangle = u(\lambda_Y(v)) = \xi(a_Y(v\otimes u)).
\end{align}
Put $R = \B_X$ and $\xi = \id_{\B_X}.$ We obtain a $\B_X$-linear map
$$
\lambda_Y : \omega(Y)\otimes \B_X \to \omega(Y)\otimes \B_X.
$$
Composing $\lambda_Y$ with the imbedding
\begin{align}\label{FactorFormula2}
\omega(Y) \to \omega(Y)\otimes \B_X,\quad v\mapsto v\otimes a_{\underline{1}}(e\otimes f),
\end{align}
where $\{e\}$ is a basis of $\omega(\underline{1})$ and
$f(e) = 1$, we provide a $\B_X$-comodule structure
$$
\rho_Y: \omega(Y) \to \omega(Y)\otimes \B_X
$$
on $\omega(Y).$
Moreover, if $Y_1,\:Y_2 \in \Ob(\Cat_X)$ and $\varphi \in \Hom(Y_1,Y_2)$ then
we have
\begin{align*}
\langle\lambda_{Y_2}\circ\omega(\varphi)(v),u\rangle &= \xi(a_{Y_2}(\omega(\varphi)(v)\otimes u)) =\xi(a_{Y_1}(v\otimes\omega(\varphi)^*(u))) =\\ &=\langle\lambda_{Y_1}(v),\omega(\varphi)^*(u)\rangle = \langle\omega(\varphi)\circ\lambda_{Y_1}(v),u\rangle
\end{align*}
for all $v \in \omega(Y_1)$ and $u \in \omega(Y_2)^*.$
Hence,
$$
\lambda_{Y_2}\circ\omega(\varphi) = \omega(\varphi)\circ\lambda_{Y_1}.
$$
Therefore, the diagram
\begin{equation*}
\begin{CD}
\omega(Y_1) @>{\rho_{Y_1}}>> \omega(Y_1)\otimes \B_X\\
@VV{\varphi}V @VV{\varphi\otimes\id_{\B_X}}V\\
\omega(Y_2)@>{\rho_{Y_2}}>> \omega(Y_2)\otimes \B_X
\end{CD}
\end{equation*}
is commutative.
This implies that
we have defined a functor
$\Cat_X \to \CoDiffT_{\B_X}$.
The forgetful functor $\CoDiffT_{\B_X} \to \Seq$ closes the commutative diagram.
\end{proof}

\begin{proposition}\label{EquivToMod} Let $\omega : \Cat \to \Seq$. Then $\omega$ defines
an equivalence of categories $\Cat_X \to \CoDiffT_{\B_X}$ taking
$\omega|_{\Cat_X}$ to the forgetful functor.
\end{proposition}
\begin{proof}
For an object $Y$ of $\Cat_X$ it is shown in Lemma~\ref{FactorsThrough}
how to put a $\B_X$-comodule structure on $\omega(Y).$
We first demonstrate that the induced functor
$$
\omega|_{\Cat_X} : \Cat_X \to \CoDiffT_{\B_X}
$$
is {\it essentially surjective}.
For a differential comodule $V$ over $\B_X$ consider the
object $\Px_X\otimes V.$ It has a subobject $S_X$ generated
by the kernel of the morphism
$$
(\varepsilon\circ\omega)\otimes\id_V : \Px_X\otimes V \to \k0\otimes
V \cong V.
$$
We then have
$$
\omega\left((\Px_X\otimes V)\big/S_X\right) = (\omega(\Px_X)\otimes
V)\big/\omega(S_X) = (\B_X\otimes V)\big/\omega(S_X) \cong V.
$$
Denote $\Px_X\otimes V\big/ S_X = V_X.$ Now, $V_X$ consists of the
collection $(Y_\alpha)$ of objects $Y_\alpha$ in $\Cat_X$ with
isomorphisms $\phi_{\beta,\alpha} : Y_\alpha \to Y_\beta$. For any
$\alpha$ we then have $\omega(Y_\alpha) \cong V$. Hence,
$\omega|_{\Cat_X}$ is essentially surjective.

We show now that $\omega|_{\Cat_X}$ is {\it full}.
For $X_1, X_2 \in \Ob(\Cat_X)$ consider the
comodules $V_1$ and $V_2$ corresponding to
$X_1$ and $X_2,$ respectively.
Let $\varphi \in \Hom(V_1,V_2).$
There is a map (see Lemma~\ref{FactorsThrough})
$$
\psi_{\varphi} := \id\otimes\varphi : \Px_X\otimes V_1 \to
\Px_X\otimes V_2,
$$
which is a morphism due to Lemmas~\ref{FromVtoHom} and \ref{HomIsom}
with $V = V_1,$ $T = \Px_X\otimes V_2,$ and $X = \Px_X.$ This induces
a morphism
$$
(V_1)_X \to (V_2)_X.
$$
By the definition of $\Sx_X$ we have surjections
$$
\id\otimes\ev : \Px_X\otimes \omega(X_i) \to X_i
$$
that induce isomorphisms
$$
(V_i)_X \to X_i.
$$
Again, $\Px_X\otimes V_i$ are the collections $(Y_\alpha^i)$ and,
hence, this $\psi$ induces a map $Y_\alpha^1 \to Y_\alpha^2$. So, we
have a morphism $X_1\to X_2$ corresponding to $\varphi$.

Since $\omega$ is {\it faithful},
we have an equivalence of categories $\Cat_X \to \Comod_{\B_X}$.
Finally, the functor $\Cat_X \to \CoDiffT_{\B_X}$ respects the functor $F$, because  $\omega$ commutes with  $F$ by definition. This establishes the required equivalence of categories.
\end{proof}

\begin{proposition}\label{InjLim} For each object $X$ of the category $\Cat$ the identification
of $\Cat_X$ with $\CoDiffT_{\B_X}$ provides an equivalence of
categories $\Cat \to \CoDiffT_A$ taking
$\omega$ to the forgetful functor from $\CoDiffT_{A}$ to $\Seq$.
\end{proposition}
\begin{proof}
Follows from Proposition~\ref{EquivToMod}. Indeed, for subcategories
$\Cat_1 \subset \Cat_2 \subset \Cat$ there is a restriction map
$\varphi_{1,2} : \IntEnd\left(\omega|_{\Cat_2}\right) \to \IntEnd(\omega|_{\Cat_1})$.
Hence, we have the induced map
${\varphi_{1,2}}^* : (\IntEnd(\omega|_{\Cat_1}))^* \to (\IntEnd(\omega|_{\Cat_2}))^*$.
Moreover, $(\Cat,\omega) = \varinjlim\left(\Cat_X, \omega|_{\Cat_X}\right)$
with respect to the injective system of inclusion maps.
\end{proof}

\begin{corollary}\label{LinearGroup} The group $G_X$, defined by $G_X(R) = \Hom(\B_X,R)$ for each $\partial$-$\k0$-algebra $R$,
is a linear differential algebraic group.
\end{corollary}
\begin{proof}
By Proposition~\ref{FinitelyGenerated} the differential Hopf algebra
$\B_X$ is finitely generated in the differential sense. Thus, according
to the proof of \cite[Proposition 12, page 914]{Cassidy} it is
a coordinate ring of a linear differential algebraic group.
\end{proof}

\begin{corollary}\label{EquivToComod}For the linear differential algebraic group $G_X$
defined above we have $$\Aut^{\otimes,\partial}\left(\omega|_{\Cat_X}\right) \cong G_X.$$
\end{corollary}
\begin{proof}
From Proposition~\ref{EquivToMod} we have equivalence of categories
$$
\Cat_X \cong \CoDiffT_{\B_X}
$$
with $\omega|_{\Cat_X}$ corresponding to the forgetful
functor $\omega'$ from $\CoDiffT_{\B_X}$ to $\Seq.$ Moreover, \cite[Theorem 2]{OvchRecoverGroup}
says that the linear differential algebraic group
$$
G_X \cong \Aut^{\otimes,\partial}(\omega'),
$$
which concludes the proof, as $\Aut^{\otimes,\partial}(\omega') \cong \Aut^{\otimes,\partial}\left(\omega|_{\Cat_X}\right).$
\end{proof}

\section{Main theorem}\label{MainSection}
\begin{theorem}\label{MainTheorem}
Let $(\Cat,\omega)$ be a neutral differential Tannakian category.
Then
$$
(\Cat,\omega) \cong \Rep_G
$$
for the differential group scheme
$$G = \Aut^{\otimes,\partial}(\omega),$$
which is
\begin{enumerate}
\item represented by the commutative associative differential Hopf
algebra $A = \varinjlim \B_X;$
\item a pro-linear differential algebraic group.
\end{enumerate}
\end{theorem}
\begin{proof}
The first statement is contained in Proposition~\ref{InjLim}.
Property~(1) follows from Corollary~\ref{EquivToComod} by taking
limits. According to Corollary~\ref{LinearGroup}
for each object $X$ of $\Cat$ the group $G_X$ is a linear differential
algebraic group. Property~(2) now follows.
\end{proof}

\section{Applications}\label{Applications}
We will show how one can apply the  theory developed above  to
parametric
linear differential
equations. The theory of parameterized linear differential equations
considers equations of the form
\[ \frac{\partial Y}{\partial x} = A(t,x)Y\]
where $A$ is an $n\times n$ matrix whose entries are functions of $x$
and
of a parameter $t$ and was developed in detail in
\cite{PhyllisMichael}.
  Formally one considers a  differential field $\K$ of characteristic
zero
with commuting derivations $\Delta = \{\partial_t, \partial_x\}$ and we
assume that $\k0 = \{ c\in \K \ | \ \partial_x(c) = 0\}$, the constants
with
respect to $\partial_x,$ forms a differentially closed
$\partial_t$-field\footnote{In \cite{PhyllisMichael},
larger
sets of derivations are considered but, for simplicity we shall
consider
one parametric derivation and one principal derivation. The subscripts
$t$
and $x$ are a convenient way of distinguishing these but we do not
assume
that we are dealing with functions of variables $t$ and $x$.}.  Given a
linear differential equation $\partial_xY = AY,$ where $A$ is an
$n\times n$
matrix with entries in $\K$, there exists a $\Delta$-differential field
extension $K$ of $\K$  having the same $\partial_x$-constants as $\K$ and
where $K$ is generated (as a $\Delta$-differential field) over $\K$ by
the
entries of an $n\times n$ invertible matrix $Z$ satisfying $\partial_xZ
=
AZ$ (Theorem 3.5, \cite{PhyllisMichael}).  This field is called the
{\em
parameterized Picard-Vessiot extension} associated to the equation and
is
unique up to differential $\K$-isomorphism.  The group of
$\K$-automorphisms
of $K$ commuting with the derivations in $\Delta$ is called the {\em
parameterized Picard-Vessiot group} of the equation and can be shown to
be
a linear differential algebraic group.  The goal of this section is to
define this group using the theory developed in the previous sections.

\subsection{Preliminaries}
Let $\K,\Delta, \k0$ be as above with $\k0$ differentially closed.  The field of constants of $\k0$ w.r.t. the differential operator
$\partial_t$ is denoted by $\k0^{\{\partial_t\}}$. A
$\partial_x$-$\K$-module $M$ is a finite dimensional $\K$-vector space with a linear
operator $\partial_x$ satisfying the product rule:
$$\partial_x(rm) = (\partial_xr)m+r\partial_xm$$ for all $r \in \K$ and $m \in M.$
Starting from Section~\ref{CategorySection} we will also use the notation $M^{(0)}$ for such a module $M.$

Let $\{e_1,\ldots,e_n\}$ be a $\K$-basis of $M$. Define the elements $a_{ij} \in \K$
by
$$
\partial_x e_i = -\sum_{j=1}^n a_{ji}e_j,
$$
where $1 \Le i \Le n.$ Let $u = a_1e_1 + \ldots + a_n e_n$.
Then $\partial_x u = 0$ iff
$$
\sum_{i=1}^n\partial_x(a_i)e_i - \sum_{i=1}^n\sum_{j=1}^na_ja_{ij}e_i = 0
$$ or, equivalently,
$$
\partial_x
\begin{pmatrix}
a_1\\
\vdots\\
a_n
\end{pmatrix}=
A\begin{pmatrix}
a_1\\
\vdots\\
a_n
\end{pmatrix},
$$
where
$A =
\begin{pmatrix}
a_{11}&\ldots&a_{1n}\\
\vdots&\ddots&\vdots\\
a_{n1}&\ldots&a_{nn}
\end{pmatrix}.$

We are going to define a formal object associated with
a parametric linear differential equation, e.g.,
$$
\partial_x y = \frac{t}{x}y
$$
which is going to be equivalent to the equation in some
sense and will allow us to recover the parametric differential
Galois group $G$ from getting first the category of all
finite dimensional differential rational representations of
$G$ and then taking all coordinate functions of representations.
We obtain the algebra $A := \k0\{G\}$ of differential algebraic
functions on $G$, then recover derivation $\partial$ on $A$,
comultiplication $\Delta : A \to A\otimes A,$ e.t.c.

\begin{remark} In what follows we will be using higher order derivation 
functors (of order 2 and higher)\footnote{see Definition~\ref{SeqCategory}.}. Formally, this is not the same as to
iterate the first order derivation functor, but:
\begin{itemize}
\item it produces objects of smaller dimension than the iterated derivatives (the derivation functor $F$ applied several times) and
\item captures the same differential information as we shall see in Proposition~\ref{iteratedishigherorder}.
\end{itemize}
Although it is not necessary to use higher order derivation functors,
we hope that this Tannakian approach will eventually produce algorithms
computing Galois groups of systems of linear differential equations with
parameters and smaller objects are more desirable for computation.
\end{remark}

The construction where $M^{(i)}$ is the $\partial_x$-module
with the matrix obtained by the $i$th prolongation (parallel to
the differential-difference construction of \cite{Charlotte,Charlotte2,CharlotteMichael1})
\begin{align}\label{AiMatrix}
A_i=\begin{pmatrix}
A &0&0&\ldots &0 \\
A_t &A& 0& \ldots & 0 \\
A_{tt}&2A_t&A&\ldots &0\\
\vdots&\vdots&\vdots&\ddots&\vdots\\
A_{t^i}&\binom{i}{1}A_{t^{i-1}}&\binom{i}{2} A_{t^{i-2}}&\ldots&A
\end{pmatrix}
\end{align}
of $A_0 = A =
\begin{pmatrix}
a_{11}&\ldots&a_{1n}\\
\vdots&\ddots&\vdots\\
a_{n1}&\ldots&a_{nn}
\end{pmatrix}$
does look like an adequate one. In the next section we will see how to
treat all these $M^{(i)},$ $i \Ge 0.$ It turns out
that we can proceed in two ways:
\begin{enumerate}
\item compute the $i$th prolongation $A_i$ of the matrix $A$,
\item compute the $i$th prolongation $M^{(i)}$ of the $\partial_x$-$\K$-module $M$. 
\end{enumerate}
Both approaches give the same $\partial_x$-structure as Proposition~\ref{ModuleToEquation} shows.

\subsection{The category $\MSeq$}\label{CategorySection}

We first give a coordinate free definition and then introduce
a canonical basis and use this constructive approach to develop our
theory.
\begin{definition}\label{Objects} We let
\begin{enumerate}
\item an object $M$ of $\MSeq$ be an
object of $\Seq$ together with a differential module structure
given by $\partial_x$ which commutes with $\partial$,
\item morphisms
between objects of $\MSeq$ be those which commute with the
action of $\partial_x:$
$$\Hom(M,N) =
\Hom_{\K[\partial_x]}(M,N)$$ for all objects $M$ and  $N$ of $\MSeq,$
\item subobjects, $\otimes,$ $\oplus,$ and $^*$ are as in the category of
differential modules \cite[Section 2.2]{Michael},
\item the functors $F^p$ are the same as in $\Seq$.
\end{enumerate}
\end{definition}

\begin{remark}\label{MSeqRemark} Recall that
\begin{align*}
&M^{(i)} = \K[\partial]_{\Le i}\otimes M^{(0)},
\end{align*}
where $M^{(0)}$ is another notation for $M.$
\end{remark}

\subsection{The equations-modules correspondence}
Derivation $\partial_x$ has to be defined
on the prolongations of $M$ in such a way that is commutes with the map $\partial : M\to M^{(1)}$, $v \mapsto \partial\otimes v$. 
\begin{lemma}\label{MatrixComputation} For an object $M$, with the $\partial_x$-structure given by
an $n\times n$  matrix $A$ with respect to an ordered basis $\{e_1,\ldots,e_n\}$,
the matrix $A_{i+1}$ corresponding to the action of $\partial_x$ on $M^{(i+1)}$ with respect to
the ordered
basis $$\{\partial^{i+1}e_1,\ldots,\partial^{i+1}e_n,\ldots,e_1,\ldots,e_n\}$$ is of the
form
$$\tilde{A}_{i+1}=\begin{pmatrix}
A &0&0&\ldots &0 \\
\binom{i+1}{1}A_t &A& 0& \ldots & 0 \\
\binom{i+1}{2}A_{tt}&\binom{i}{1}A_t&A&\ldots &0\\
\vdots&\vdots&\vdots&\ddots&\vdots\\
A_{t^{i+1}}& \binom{i}{i} A_{t^i}&A_{t^{i-1}}&\ldots&A
\end{pmatrix}.$$
\end{lemma}
\begin{proof}
Let $$A_0 = A =
\begin{pmatrix}
a_{11}&\ldots&a_{1n}\\
\vdots&\ddots&\vdots\\
a_{n1}&\ldots&a_{nn}
\end{pmatrix}$$
be the matrix for $\partial_x$ in the module $M^{(0)}$.
For $i=1$ the basis of the module $M^{(1)}$ is $\{\partial e_1,\ldots,\partial e_n,e_1,\ldots,e_n\}$.
For $k$, $1 \Le k \Le n$, we have
\begin{align*}
\partial_x(\partial e_k) =\partial(\partial_x(e_k)) =
\partial\left(-\sum_{p=1}^n a_{pk}e_p\right) &=
-\sum_{p=1}^n a_{pk}\partial e_p -\sum_{p=1}^n (\partial_ta_{pk})e_p,\\
\partial_x(e_k) &= -\sum_{p=1}^n a_{pk}e_p.
\end{align*}
Assume the result for $i=m$. We have
$\{\partial^{m+1}e_1,\ldots,\partial^{m+1}e_n,\ldots,e_1,\ldots,e_n\}$
as the ordered basis of $M^{(m+1)}$.
For each $l,$ $1\Le l \Le n,$ we have
\begin{align*}
\partial_x&(\partial^{m+1}e_l) = \partial(\partial_x(\partial^me_l))=\\
&= \partial\left(-\sum_{q=0}^{m}\sum_{r=1}^n\binom{m}{q}(\partial_t^q a_{rl})\partial^{m-q}e_r\right)=\\
&= -\sum_{q=0}^{m}\sum_{r=1}^n\binom{m}{q}(\partial_t^{q+1} a_{rl})\partial^{m-q}e_r
-\sum_{q=0}^{m}\sum_{r=1}^n\binom{m}{q}(\partial_t^q a_{rl})\partial(\partial^{m-q}e_r) = \\
&= -\sum_{q=0}^{m}\sum_{r=1}^n\binom{m}{q}(\partial_t^{q+1} a_{rl})\partial^{m+1-(q+1)}e_r
-\sum_{q=0}^{m}\sum_{r=1}^n\binom{m}{q}(\partial_t^q a_{rl})\partial^{m+1-q}e_r =\\
&= -\sum_{q=1}^{m+1}\sum_{r=1}^n\binom{m}{q-1}(\partial_t^q a_{rl})\partial^{m+1-q}e_r
-\sum_{q=0}^{m}\sum_{r=1}^n\binom{m}{q}(\partial_t^q a_{rl})\partial^{m+1-q}e_r =\\
&= -\sum_{q=0}^{m+1}\sum_{r=1}^n\binom{m+1}{q}(\partial_t^q a_{rl})\partial^{m+1-q}e_r,
\end{align*}
because $\binom{m}{q-1}+ \binom{m}{q} = \binom{m+1}{q}$.
\end{proof}

\begin{proposition}\label{ModuleToEquation} For an object $M$ of $\MSeq$ there exists an
ordered basis of each $M^{(i)}$ such that the matrix $A_{i}$
corresponding to the action of $\partial_x$ w.r.t. this basis
is of the form
$$A_i=\begin{pmatrix}
A &0&0&\ldots &0 \\
A_t &A& 0& \ldots & 0 \\
A_{tt}&\binom{2}{1}\cdot A_t&A&\ldots &0\\
\vdots&\vdots&\vdots&\ddots&\vdots\\
A_{t^i}&\binom{i}{1}A_{t^{i-1}}&\binom{i}{2} A_{t^{i-2}}&\ldots&A
\end{pmatrix},$$
where  $$A_0 = A =
\begin{pmatrix}
a_{11}&\ldots&a_{1n}\\
\vdots&\ddots&\vdots\\
a_{n1}&\ldots&a_{nn}
\end{pmatrix}$$
is the matrix for $\partial_x$ in the module $M$.
\end{proposition}
\begin{proof} According to Lemma~\ref{MatrixComputation} for each
$M^{(i)}$ there exists an ordered basis $$\{\partial^ie_1,\ldots,\partial^ie_n,\ldots,e_1,\ldots, e_n\}$$ such that the matrix corresponding
to $\partial_x$ is
$$\tilde{A}_i=\begin{pmatrix}
A &0&0&\ldots &0 \\
\binom{i}{1}\cdot A_t &A& 0& \ldots & 0 \\
\binom{i}{2}\cdot A_{tt}&\binom{i-1}{1}\cdot A_t&A&\ldots &0\\
\vdots&\vdots&\vdots&\ddots&\vdots\\
A_{t^i}&\ldots&\ldots&\ldots&A
\end{pmatrix}.$$
Denote this basis by $\{f_1,\ldots,f_{(i+1)\cdot n}\}$. Note that if
we change a basis:
$$\begin{pmatrix}
g_1\\
\vdots\\
g_{(i+1)\cdot n}
\end{pmatrix}=
C\cdot\begin{pmatrix}
f_1\\
\vdots\\
f_{(i+1)\cdot n}
\end{pmatrix}$$
and all entries of the matrix $C$ are constants then the corresponding
matrix $A_i$ must be replaced by $(C^T)^{-1}A_iC^T$. Consider the basis $\{g_1,\ldots,g_{(i+1)\cdot n}\}$ given by the following change-of-basis matrix
$$
C := \begin{pmatrix}
E &0&0&\ldots &0 \\
E &\binom{i}{1}\cdot E& 0& \ldots & 0 \\
E &\binom{i-1}{1}\cdot E&\binom{i}{2}\cdot E&\ldots &0\\
\vdots&\vdots&\vdots&\ddots&\vdots\\
E &2\cdot E&3\cdot E&\ldots &0\\
0&E&E&\ldots&E
\end{pmatrix}^T,
$$
where $E$ is the $n\times n$ identity matrix. Indeed, for $l,$ $1 \Le l \Le n$ we have
\begin{align*}
\partial_x g_l &= \partial_x (f_l+f_{l+n}+\ldots+f_{l+(i-1)\cdot n}) =\\
&= \partial_x(\partial^{i-1}e_l+\partial^{i-2}e_l+\ldots+e_l) = \\
&= \sum_{p=0}^{i-1}\left(-\sum_{q=p}^i\sum_{k=1}^n\binom{i-p}{q-p}(\partial_t^{q-p}a_{kl})\partial^{i-1-q}e_k\right) =\\
&=- \sum_{p=0}^i\sum_{k=1}^n\left(\sum_{q=0}^p\binom{i}{q}(\partial_t^qa_{kl})\right)\partial^{i-1-p}e_k =\\
&= -\sum_{p=0}^i\sum_{k=1}^n (\partial_t^pa_{kl})\left(\sum_{q=p}^i\binom{i-(q-p)}{p}\partial^{i-1-q}e_k\right)=\\
&= -\sum_{p=0}^i\sum_{k=1}^n (\partial_t^pa_{kl})g_{k+p\cdot n}
\end{align*} and for all $k,$ $1 \Le k \Le n,$ we have
\begin{align*}
\partial_x&(g_{k+p\cdot n}) = \partial_x\left(\sum_{q=p}^i\binom{i-(q-p)}{p}f_{k+q\cdot n}\right)=\\
&=\partial_x\left(\sum_{q=p}^i\binom{i-(q-p)}{p}\partial^{i-1-q}e_k\right)=\\
&= \sum_{q=p}^i\binom{i-(q-p)}{p}\partial_x \partial^{i-1-q}e_k =\\
&=-\sum_{q=p}^i\binom{i-(q-p)}{p}\sum_{l=1}^n
(\partial_t^{q-p}a_{lk})\partial^{i-1-q}e_l =\\
&=-\sum_{r=0}^{i-p}\sum_{l=1}^n\binom{p+r}{r}
(\partial_t^r a_{lk})\left(\sum_{s=p+r}^i
\frac{\binom{i-s+p+r}{p}\binom{i-p+p+r-s}{r}}
{\binom{p+r}{r}}\partial^{i-1-s}e_l\right)=\\
&=-\sum_{r=0}^{i-p}\binom{p+r}{r}
\sum_{l=1}^n(\partial_t^r a_{lk})\left(\sum_{s=p+r}^i
\binom{i-s+p+r}{p+r}\partial^{i-1-s}e_l\right) =\\
&=-\sum_{r=0}^{i-p}\binom{p+r}{r}
\sum_{l=1}^n(\partial_t^r a_{lk})g_{l+(p+r)\cdot n},
\end{align*}
because
\begin{align*}
\frac{\binom{i-s+p+r}{p}\binom{i-p+p+r-s}{r}}{\binom{p+r}{r}}
&= \frac{(i-s+p+r)!(i+r-s)!r!p!}{p!(i-s+r)!r!(i-s)!(p+r)!}=\\
&=\frac{(i-s+p+r)!}{(i-s)!(p+r)!}=\binom{i-s+p+r}{p+r}.
\end{align*}
\end{proof}

\subsection{Covariant solution space} Fix an object $M$
of $\MSeq$.
Let $K$ be a parametric Picard-Vessiot extension of $\K$ for the equation
$$
\partial_x Y = AY,
$$
where $A$ is the matrix corresponding
to the $\K$-finite dimensional $\K[\partial_x]$-module $M.$ Then the
covariant solution space is $$V = \ker(\partial_x,K\otimes M),$$
which is a vector space over the field $\k0$.
We let $$V^{(i)} = \ker\left(\partial_x,K\otimes M^{(i)}\right),$$ which is an
$(i+1)\cdot\dim M^{(0)} = (i+1)\cdot n$-dimensional  vector space over $\k0$.
Let $Y \in \GL_n(K)$ be a solution matrix for $V^{(0)}$. Then,
the columns of the matrix
$$
Y_i := \begin{pmatrix}
Y &0&0&\ldots &0 \\
Y_t & Y& 0& \ldots & 0 \\
Y_{tt} & Y_t& Y&\ldots &0\\
\vdots&\vdots&\vdots&\ddots&\vdots\\
Y_{t^{i-1}} &Y_{t^{i-2}}&Y_{t^{i-3}}&\ldots &0\\
Y_{t^i}&Y_{t^{i-1}}&Y_{t^{i-2}}&\ldots&Y
\end{pmatrix}
$$
form a $\k0$-basis for $V^{(i)},$ since $Y_i \in \GL_{n\cdot(i+1)}(K).$
Hence,
for each $M$ we have a map sending $M^{(i)} \mapsto V^{(i)}$ and
$V^{(i)}$ can be identified with the vector space $\k0[\partial]_{\Le i}\otimes V^{(0)}$ by Proposition~\ref{ModuleToEquation}. Moreover, \eqref{difftensor}, \eqref{diffexact}, \eqref{diffdual}, and \eqref{diffdualomega} are satisfied. So, there
is a functor $\omega : M \mapsto V$ from the category $\MSeq$ to the category $\Seq$.

\subsection{Some examples}

\begin{example} Consider a matrix $A \in \Mn_n\left(\k0^{\{\partial_t\}}\right)$
and the differential equation $$\partial_x Y = AY.$$
Then, $V_i$ is the $\partial_x$-differential module with
the matrix
$$A_i = \diag(A,\ldots,A)
$$
and the $(i+1)$th solution matrix is
$$
Y_i = \diag (Y,\ldots,Y).
$$
So, we do not bring anything new to the usual differential Galois
theory.
\end{example}

\begin{example} Consider the differential equation
$$
\partial_x y = \frac{t}{x}y.
$$
Then,
\begin{align*}
A_0 &= \begin{pmatrix}\frac{t}{x}\end{pmatrix},\\
A_1 &=
\begin{pmatrix}
\frac{t}{x} & 0\\
\frac{1}{x} & \frac{t}{x}
\end{pmatrix},\\
&\vdots\\
A_i &=
\begin{pmatrix}
\frac{t}{x} & 0 & 0 & 0&\ldots & 0\\
\frac{1}{x} & \frac{t}{x} & 0 & 0&\ldots &0\\
0 & \frac{2}{x} & \frac{t}{x} & 0 & \ldots &0\\
\vdots&\vdots&\vdots&\vdots&\ddots&\vdots\\
0&0&0&0&\ldots&\frac{t}{x}
\end{pmatrix},\\
&\vdots
\end{align*}
and the solution matrices are
\begin{align*}
Y_0 &= \begin{pmatrix}x^t\end{pmatrix},\\
Y_1 &=
\begin{pmatrix}
x^t & 0\\
x^t\log x & x^t
\end{pmatrix},\\
&\vdots\\
Y_i &=
\begin{pmatrix}
x^t & 0 & 0 & 0&\ldots & 0\\
x^t\log x & x^t & 0 & 0&\ldots &0\\
x^t(\log x)^2 & x^t\log x & x^t & 0 & \ldots &0\\
\vdots&\vdots&\vdots&\vdots&\ddots&\vdots\\
x^t(\log x)^i&x^t(\log x)^{i-1}&x^t(\log x)^{i-2}&x^t(\log x)^{i-3}&\ldots&x^t
\end{pmatrix},\\
&\vdots
\end{align*}

\end{example}

\subsection{Fiber functor}\label{FibreFunctor}

\begin{proposition} The functor $\omega$ is an exact $\k0$-linear
faithful tensor covariant functor intertwining with the functors $F^p$ and preserving the map $\partial$.
\end{proposition}
\begin{proof}
We have:
\begin{itemize}
\item Consider a  short exact sequence
\begin{equation*}
\begin{CD}
0 @>{}>> L @>{\varphi_1}>>M @>{\varphi_2}>>N @>{}>>0\\
\end{CD}
\end{equation*}
Since $\K$ is a field,
the sequence
\begin{equation}\label{first}
\begin{CD}
0 @>{}>> K\otimes_k L @>{\id_K\otimes\varphi_1}>>K\otimes_k M @>{\id_K\otimes\varphi_2}>>K\otimes_k N @>{}>>0\\
\end{CD}
\end{equation}
is also exact. We have
\begin{align*}
\dim_K K\otimes M-\dim_K K\otimes L &= \dim_K K\otimes N,\\
\dim_{\k0}\ker\left(\partial_x,K\otimes M\right) &= \dim_K K\otimes M,\\
\dim_{\k0}\ker\left(\partial_x,K\otimes L\right) &= \dim_K K\otimes L,\\
\dim_{\k0}\ker\left(\partial_x,K\otimes N\right) &= \dim_K K\otimes N.
\end{align*}
Consider the following short sequence:
\begin{equation}\label{second}
\begin{CD}
0 @>{}>> \ker\left(\partial_x,K\otimes_{\K} L\right) @>{}>>\ker\left(\partial_x,K\otimes_{\K} M\right)
@>{}>>\\
@>{}>>\ker\left(\partial_x,K\otimes_{\K} N\right) @>{}>>0\\
\end{CD}
\end{equation}

The first and second elements of this sequence are subsets of
the corresponding elements of the sequence~\eqref{first} and, hence, sequence~\eqref{second} is exact at those elements. Since
$\ker(\partial_x,K\otimes M)/\ker(\partial_x,K\otimes L)$ is
isomorphic to a subspace of $\ker(\partial_x,K\otimes N)$ and
according to the dimensional equalities we have
$$
\dim_{\k0}\ker\left(\partial_x,K\otimes M\right)/\ker\left(\partial_x,K\otimes L\right) =
\dim_{\k0}\ker\left(\partial_x,K\otimes N\right),
$$
the sequence~\eqref{second} is exact at the third element.
\item We have $\varphi = a\varphi_1 + b\varphi_2$ commutes with $\partial_x$ for all $\varphi_i \in \Hom(M,N)$
and $a,b \in \k0.$ Hence, $\varphi$ is mapped to $a\cdot\omega(\varphi_1)+b\cdot \omega(\varphi_2).$ Thus, $\omega$ is
$\k0$-linear.
\item A non-zero module provides a non-zero solution.
\item For modules  $M$ and $N$
we only need to show that
$$
\ker\left(\partial_x,K\otimes M\right)\otimes\ker\left(\partial_x,K\otimes N\right)=
\ker\left(\partial_x,K\otimes M\otimes N\right),
$$
but this follows as in the usual differential case by
taking $K$-bases $\{e_i\}$ and $\{f_j\}$ of $K\otimes M$ and
$K\otimes N$, respectively, with $\partial_x(e_i) = \partial_x(f_j) =0$.
\item Let $\varphi \in \Hom(M,N)$. Then it is mapped to
the morphism $\id\otimes\varphi$ of the spaces $\ker(\partial_x,K\otimes M)$
and $\ker(\partial_x,K\otimes N)$.
\end{itemize}
\end{proof}

\begin{corollary}\label{MIsTannaka} The category $\MSeq$ together with the functor $\omega$
is a neutral differential Tannakian category (Definition~\ref{NeutralParametricTannakian}).
\end{corollary}

For $M \in \Ob(\MSeq)$ denote by $\Cat_M$ the category generated by the module $M$ and all its derivatives $M^{(i)}$ and the dual $M^*$ using the operations $\oplus$, $\otimes$, and
subquotient.  We say that $\Cat_M$ is {\it differentially generated} by $M$.

\begin{corollary}\label{RecoverDiffGalois}For a parametric linear differential equation $$\partial_xY=A(x,t)Y$$
a parametric differential Galois  group $G$ is the functor $\Aut^{\otimes,\partial}(\omega)$,
where $\omega : \Cat_M \to \Vect_{\k0}$ is the functor associated with $M\in \MSeq$ and
the action of $\partial_x$ defined by the matrix $A \in \Mn_n(\K).$
\end{corollary}
\begin{proof}
From Theorem~\ref{MainTheorem} and Corollary~\ref{MIsTannaka} it follows
that the category $\Cat_M$ together with the functor $\omega$ form a neutral differential Tannakian
category.
Moreover, $\Cat_M$ is equivalent to the category of representations
of the linear differential algebraic group $\Aut^{\otimes,\partial}(\omega)$. It remains
to show that for the parametric Galois group $G$ its category of representations
is equivalent to $\Cat_M$.

Similar to \cite[Theorem 2.33]{Michael}, in order to get this equivalence
we notice the following. Consider a parameterized Picard-Vessiot
extension $K$ of $\K$  corresponding to the equation $\partial_x Y = AY$. This field $K$ is generated over $\K$ by
a fundamental set of solutions of the differential equation and $G$ acts on the
$\k0$-linear space generated by these solutions. So, $G$ acts on $K$
and commutes with $\partial_x$. Hence, if for $f\otimes m \in K\otimes M$
and $\sigma \in G$, which is an automorphism of $K$, we let
$$
\sigma (f\otimes m) = \sigma(f)\otimes m,
$$
then the linear differential algebraic group $G$ acts faithfully on $\omega(M)$,
which is an object of the category $\Seq$. Indeed, recall that if one fixes a $\K$-basis
of $M$, one gets a fundamental solution matrix, on which $\sigma$ acts  by multiplication by an invertible matrix with coefficients in $\k0$. 

From \cite[Proposition 3]{OvchRecoverGroup} it follows that the category
$\Rep_G$ is {\it differentially} generated by a faithful representation $V$
of $G$. We take $V = \omega(M)$ and, hence, for every representation $W$ there
exists an object $M_W$ in $\Cat_M$ such that $\omega(M_W) = W$. Therefore, the
functor $\omega$ is essentially surjective. We need to show that
$$
\Hom(M,N) \to \Hom(S(M),S(N))
$$
is a bijection.
We have:
\begin{align*}
\Hom_{\K[\partial_x]}\left(M,N\right) = \ker\left(\partial_x,M^*\otimes N\right) =\Hom_{\K[\partial_x]}\left(\k0,M^*\otimes N\right).
\end{align*}
So, it is enough to show the bijection for $M = \K.$ We have
$\omega(\K) = \k0$, which is the trivial representation
of the group $G.$ Moreover,
$$
\Hom_{\K[\partial_x]}\left(\k0,N\right) = \left\{n \in N\:\big|\:\partial_x n = 0\right\}
$$
and
$$
\Hom_{\k0}(\k0,\omega(N)) = \left\{v \in \ker\left(\partial_x,K\otimes N\right)\:\big|\: gv=v\ \text{for all}\ g \in G\right\}.
$$
But the fixed points $K^G$ are $\K$ (see \cite[Theorem 3.5]{PhyllisMichael}). Hence, if $g(f\otimes n) = f\otimes n$
for $f \in K$ and $n \in N$ for all $g\in G$ then $f \in \K$, and
$$
\left(K\otimes N\right)^G = N.
$$
Thus, we have shown a bijection between
$\Hom(M,N)$ and $\Hom(\omega(M),\omega(N))$.
\end{proof}

We will finally show that our more compact and, therefore, computationally friendlier
way of recovering $G$  is equivalent
to applying the derivation functor $F$ iteratively. Recall that if one applies $F$ twice,
one multiplies the dimension by $4$. In our construction we multiply it only by $3$.

\begin{proposition}\label{iteratedishigherorder} 
For $M \in \MSeq$ the category $(\Cat_M, \omega)$ coincides with the neutral differential
Tannakian category $(\Cat,\omega)$ generated by $M$ using $\oplus$, $\otimes$, $*$, and
iterated application of the first order derivation functor $F$.
\end{proposition}
\begin{proof} 
First note that the category $\Cat_M$ is a subcategory of $\Cat$. Indeed,
denote the $i$-th iteration of $F$ on $M$ by $M^{[i]}$. Then, $M^{(i)}$
can be embedded into $M^{[i]}$. For $i=1$ we have $M^{[1]}=M^{(1)}$. For simplicity, we will show an embedding $M^{(i)} \to M^{[i]}$ for $i=2$. The general case can be shown by induction on $i$. Let $m \in M$. We will
then define an embedding $E$ by
\begin{align*}
1\otimes  m &\mapsto 1\otimes 1\otimes m,\\
\partial\otimes  m &\mapsto \dfrac{\partial\otimes 1\otimes m+1\otimes\partial\otimes m}{2},\\
\partial^2\otimes m &\mapsto \partial\otimes\partial\otimes m.
\end{align*}
This map is $\K$-linear since it preserves the degree of $\partial$. We will show that
it commutes with $\partial_x$. Suppose a  basis of $M$ is fixed and $A\in\Mn_n(\K)$
defines the action of $\partial_x$ with respect to this basis. Lemma~\ref{MatrixComputation} provides the matrices
$$
A^{(2)} :=\begin{pmatrix}
A&0&0\\
2A_t&A&0\\
A_{tt}&A_t&A
\end{pmatrix}\quad \text{and}
\quad A^{[2]} := 
\begin{pmatrix}
A&0&0&0\\
A_t&A&0&0\\
A_t&0&A&0\\
A_{tt}&A_t&A_t&A
\end{pmatrix}
$$
corresponding to the prolongations of $\partial_x$ to $M^{(2)}$ and $M^{[2]}$, respectively, from where the commutation of $\partial_x$ and $E$ follows. Thus,
$E \in \Hom\left(M^{(2)},M^{[2]}\right)$.

By Corollary~\ref{RecoverDiffGalois}, the category $(\Cat,\omega)$ is a neutral
differential Tannakian category and is, therefore, equivalent to $\Rep_G$, where
$G$ is the parametrized differential Galois group for the module $M$. By \cite[Proposition 2]{OvchRecoverGroup}, the image of $\Cat_M$ under this equivalence
coincides with the whole $\Rep_G$, and, thus, $\Cat_M = \Cat$.
\end{proof}

\section{Conclusions} Starting with a parametric linear differential equation
one first constructs the parametric differential module $M$ associated with it.
To apply  linear algebra operations 
and the prolongation functor to it is the same as to do them
with its solution space. In this way we construct the category $\Cat_M$ differentially generated by $M$
with a fiber functor $S : \Cat_M \to \Seq$. From this data one recovers
all finite dimensional differential representations and, hence, the parametric differential Galois group
of the equation. 

\section{Acknowledgements}
The author is highly grateful to his advisor Michael Singer,
to Bojko Bakalov, and Daniel Bertrand
for extremely helpful comments and support. Also, the author
thanks Pierre Deligne, Christian Haesemeyer, Moshe Kamensky, Claudine Mitschi,  Jacques Sauloy, the participants of Kolchin's Seminar in New York,
Sergey Gorchinsky, and the referees for their important suggestions.

\end{document}